\documentclass[11pt,twoside]{article}
\usepackage{txfonts}
\usepackage{amssymb}
\usepackage{graphicx}
\usepackage{cases}
\usepackage{indentfirst}
\newtheorem{thm}{Theorem}[section]
\newtheorem{lem}[thm]{Lemma}
\newtheorem{cor}[thm]{Corollary}
\newtheorem{conj}[thm]{Conjecture}

\def\qed{\hfill \rule{2.5mm}{2.5mm}}
\date{}

\begin{document}

\title{\LARGE{List version of ($p$,1)-total labellings} \thanks{This work is supported by IIFSDU, NNSF(10871119) and
RSDP(200804220001) of China.}}

\author{\Large{Yong Yu, Guanghui Wang, Guizhen Liu}
\thanks {Corresponding author. E-mail address: gzliu@sdu.edu.cn.}\\ [5pt] \small{School of Mathematics,
Shandong University, Jinan 250100, P.R.China}}
\maketitle

\def\abstractname{Abstract}
\begin{abstract}
The ($p$,1)-total number $\lambda_p^T(G)$ of a graph $G$ is the
width of the smallest range of integers that suffices to label the
vertices and the edges of $G$ such that no two adjacent vertices
have the same label, no two incident edges have the same label and
the difference between the labels of a vertex and its incident edges
is at least $p$. In this paper we consider the list version. Let
$L(x)$ be a list of possible colors for all $x\in V(G)\cup E(G)$.
Define $C_{p,1}^T(G)$ to be the smallest integer $k$ such that for
every list assignment with $|L(x)|=k$ for all $x\in V(G)\cup E(G)$,
$G$ has a ($p$,1)-total labelling $c$ such that $c(x)\in L(x)$ for
all $x\in V(G)\cup E(G)$. We call $C_{p,1}^T(G)$ the ($p$,1)-total
labelling choosability and $G$ is list $L$-($p$,1)-total labelable.
\par In this paper, we present a conjecture on the upper bound of
$C_{p,1}^T$. Furthermore, we study this parameter for paths and
trees in Section 2. We also prove that $C_{p,1}^T(K_{1,n})\leq
n+2p-1$ for star $K_{1,n}$ with $p\geq2, n\geq3$ in Section 3 and
$C_{p,1}^T(G)\leq \Delta+2p-1$ for outerplanar graph with
$\Delta\geq p+3$ in Section 4.

\vspace{3mm}

\end{abstract}



\vspace{5mm}

\bigskip

{\section{\bf Introduction}}
\bigskip
 In this paper, the term \emph{graph} is used to denote a simple
connected graph $G$ with a finite vertex set $V(G)$ and a finite
edge set $E(G)$. The degree of a vertex $v$ in $G$ is the number of
edges incident with $v$ and denoted by $d_G(v)$. We write
$\delta(G)=\min\{d_G(v): v\in V(G)\}$ and $\Delta(G)=\max\{d_G(v):
v\in V(G)\}$ to denote the minimum degree and maximum degree of $G$,
respectively. We sometimes write $V, E, d(v), \Delta, \delta$
instead of $V(G), E(G), d_G(v), \Delta(G), \delta(G)$, respectively.
A function $L$ is called an \emph{assignment} for a graph $G$ if it
assigns a list $L(x)$ of possible labels (or colors) to each element
$x\in V(G)\cup E(G)$. A \emph{$k$-assignment} is a list assignment
where all lists have the same cardinality $k$, that is, $|L(x)|=k$
for all $x\in V(G)\cup E(G)$. We shall assume throughout that the
labels (or colors) are natural numbers.  Our terminology and
notation will be standard except where indicated. Readers are
referred to \cite{2} for undefined terms.
\par Let $p$ be a  nonnegative integer. A \emph{$k$-($p$,1)-total
labelling} of a graph $G$ is a function $c$ from $V(G)\cup E(G)$ to
the color set $\{0,1,\cdots,k\}$ such that $c(u)\neq c(v)$ if $uv\in
E(G)$, $c(e)\neq c(e')$ if $e$ and $e'$ are two adjacent edges, and
$|c(u)-c(e)|\geq p$ if vertex $u$ is incident to the edge $e$. The
minimum $k$ such that $G$ has a $k$-($p$,1)-total labelling is
called the \emph{($p$,1)-total labelling number} and denoted by
$\lambda_p^T(G)$. Let us denote by $\chi_{p,1}^T(G)$ the minimum
number of colors(labels) needed for an ordinary ($p$,1)-total
labelling for describing conveniently in this paper. Obviously, we
have $\chi_{p,1}^T(G)=\lambda_p^T(G)+1$. When $p=1$, the (1,1)-total
labelling is the well-known total coloring of graphs, and
$\chi_{1,1}^T(G)=\chi''(G)$ where $\chi''(G)$ denotes the total
chromatic number.
\par Here we present the concept list ($p$,1)-total labelling.
Suppose $L$ is an assignment for a graph $G$. If $G$ has a
($p$,1)-total labelling $c$ such that $c(x)\in L(x)$ for all $x\in
V(G)\cup E(G)$, then we say that $c$ is an $L$-($p$,1)-total
labelling of $G$, and $G$ is $L$-($p$,1)-total labelable.
Furthermore, if $G$ is $L$-($p$,1)-total labelable for any $L$ with
$|L(x)|=k$ for each $x\in V(G)\cup E(G)$, we say that $G$ is
$k$-($p$,1)-total choosable. The ($p$,1)-total labelling
choosability, denoted by $C_{p,1}^T(G)$, is the minimum $k$ such
that $G$ is $k$-($p$,1)-total choosable.
\par Obviously, this concept is a common generalization of list
colorings and ($p$,1)-total labellings. The ($p$,1)-total labelling
of graphs was introduced by Havet and Yu \cite{6}. It was shown that
$\lambda_p^T(G)\leq2\Delta+p-1$ for any graph $G$, and if
$\Delta\geq3$ then $\lambda_p^T(G)\leq2\Delta$, if $\Delta\geq5$ is
odd then $\lambda_p^T(G)\leq2\Delta-1$. The special cases for $p=2$
were also investigated in this paper. Some kind of special graphs
have also been studied, e.g., complete bipartite graphs for $p=2$
\cite{9}, planar graphs \cite{1}, trees for $p=2$ \cite{12}, graphs
with a given maximum average degree \cite{10}, complete graphs
\cite{6}, etc. In [6], Havet and Yu gave a conjecture that
$\lambda_p^T(G)\leq\Delta+2p-1$ for any graph $G$, which extends the
well known Total Coloring Conjecture in which $p=1$.
\par The \emph{incidence graph} of a graph $G$, denoted by $S_I(G)$,
is the graph obtained from $G$ by replaced each edge by a path of
length 2. Motivated by the Frequency Channel Assignment Problem,
Griggs and Yeh \cite{5} first introduced the $L(2,1)$-labelling of
graphs. This notion was subsequently extended to a general form,
named as $L(p,q)$-labelling of graphs. The $L(p,q)$-labelling,
especially the $L$(2,1)-labelling, of graphs have been studied
rather extensively in recent years. Kohl et al. \cite{4}
investigated the list version of $L(p,q)$-labellings and obtained
some interesting results. As mentioned in \cite{6}, the
$L(p,1)$-labelling of $S_I(G)$ is equivalent to the ($p$,1)-total
labelling of graph $G$. We still noticed that the ($p$,1)-total
labelling is a special case of an $[r, s, t]$-coloring of graphs
with $r=s=1,t=p$, which was introduced in \cite{7}. Hence it is easy
to
see :\\

{\bfseries Observation 1.} Let $G$ be a graph. Then
$$\chi_l^{p,1}(S_I(G))=C_{p,1}^T(G)=\chi_l^{1,1,p}(G),$$ where
$\chi_l^{p,1}(G)$ and $\chi_l^{1,1,p}(G)$ denote the minimum number
$k$ such that $G$ is $k$-$L(p,1)$-labelling choosable and $k$-$[1,
1, p]$-coloring choosable, respectively.
\par In Section 2, we give some general bounds for $C_{p,1}^T$ for paths and
trees. After that, we present a conjecture on the upper bound of
$C_{p,1}^T$ for any graph $G$.
\par In Section 3, we show that $C_{p,1}^T(K_{1,n})\leq n+2p-1$ where $p\geq2,
n\geq3$.
\par In Section 4, we discuss the value for $C_p^T$ for outerplanar
graphs. We prove that $C_{p,1}^T(G)\leq \Delta+2p-1$ for all
outerplanar graph $G$ with $\Delta\geq p+3$, and we conjecture that
the upper bound is still true without the maximum degree
restriction.\\

\bigskip

{\section{\bf Basic results on $C_{p,1}^T$}}
\par At first, by using Observation 1 we try to give some bounds for paths and
trees. Then we give a conjecture on the upper bound for any graph
$G$.
\medskip
\begin{lem}
{\upshape ([8]or[4])} Let $P_k$ be a path with $k$ vertices. Then
$$\chi^{p,1}(P_k)=\left\{
\begin{array}{ll}
p+1, & \hbox{$k=2$;} \\
p+2, & \hbox{$k=3,4$;} \\
p+3, & \hbox{$k\geq5$.}
\end{array}
\right.$$
\end{lem}

\begin{thm}
Let $P_k$ be a path with $k$ vertices. Then
$$\chi_{p,1}^T(P_k)=\left\{
\begin{array}{ll}
p+2, & \hbox{$k=2$;} \\
p+3, & \hbox{$k\geq3$.}
\end{array}
\right.$$
\end{thm}

\noindent{\bf Proof.} Let $S_I(P_k)$ be the incidence graph of
$P_k$, then $S_I(P_k)=P_{k'}$ is still a path with $k'=2k-1$. By
Observation 1 and Lemma 2.1, we have
$\chi_{p,1}^T(P_2)=\chi^{p,1}(P_3)=p+2$, and when $k\geq3$ we have
$\chi_{p,1}^T(P_k)=\chi^{p,1}(P_{k'})=p+3$ since $k'=2k-1\geq5$.\qed

\begin{lem} {\upshape ([8])} Let $P_k=v_1\cdots v_k$ be a path and $k>2p$. Then
we have $2p\leq\chi^{p,1}_l(P_{k})\leq2p+1$.
\end{lem}

\begin{thm}
Let $P_k=v_1\cdots v_k$ be a path. Then $C_{p,1}^T(P_k)\leq2p+1$.
Moreover, if $k>p$, then we have $2p\leq C_{p,1}^T(P_k)\leq2p+1$.
\end{thm}

\noindent{\bf Proof.} $C_{p,1}^T(P_k)\leq2p+1$ is obvious since we
can color the vertices and edges of the path sequentially in its
order by a greedy algorithm. When $k>p$, an analogous argument with
the proof in Theorem 2.2 shows that
$C_{p,1}^T(P_k)=\chi_l^{p,1}(S_I(P_k))$. Then by Lemma 2.3 we have
$2p\leq C_{p,1}^T(P_k)\leq2p+1$.\qed

\vspace{5mm}

\par When $p=2$, we have $C_{2,1}^T(P_k)\leq5$ with $k\leq3$ by
Theorem 2.4. By the definition of $C_{p,1}^T$, it is easy to see
that $C_{p,1}^T(G)\geq \chi_l^{p,1}(G)$. Then $C_{2,1}^T(P_k)\geq
\chi_l^{2,1}(P_k)=5$ by Theorem 2.2. Therefore,
$C_{2,1}^T(P_k)=5=2p+1$ when $k\geq3$. So the upper bound of Theorem
2.4 is tight.

\begin{lem} {\upshape ([8])} For all trees $T$, all $d$ and all $s\geq1$, we
have $\chi_l^{d,s}(T)\leq 2d-1+s\Delta$.
\end{lem}

\begin{thm} Let $T_n$ be a tree with $n$ vertices. Then we have $C_{p,1}^T(T_n)\leq
\Delta+2p-1$.
\end{thm}

\noindent{\bf Proof.} Let $S_I(T_n)$ be the incidence graph of
$T_n$. $S_I(T_n)$ is still a tree with $n'=2n-1$ vertices and
$\Delta(T_{n'})=\Delta(T_n)$. By Lemma 2.5, let $d=p,s=1$ we have
$\chi_l^{p,1}(T_{n'})\leq \Delta+2p-1$. Therefore, by Observation 1
we obtain $C_{p,1}^T(T_n)=\chi_l^{p,1}(T_{n'})\leq \Delta+2p-1$.\qed

\begin{lem} {\upshape ([8])} If $T$ is a tree with maximum degree $\Delta,
p\leq\Delta$, and there is a vertex $v\in V(G)$ such that $v$ and
all of its neighbors have degree $\Delta$, then
$\chi_l^{p,1}(T_{n})=\Delta+2p-1$.
\end{lem}

\par By Lemma 2.7, if $T=P_n, n\geq5$ and $p=2$,
$C_{2,1}^T(T_n)=\chi_l^{2,1}(S_I(T_n))=\chi_l^{p,1}(T_{n'})=\Delta(T_{n'})+2p-1=\Delta(T_n)+2p-1$.
That is to say, the upper bound of Theorem 2.6 is also tight.
\par It is known to all that for list version of edge colorings and
total colorings there are list edge coloring conjecture (LECC) and
list total coloring conjecture (LTCC) as follows:\\
$$(1) \chi_l'(G)=\chi'(G);\  (2) \chi_l''(G)=\chi''(G).$$
Therefore, it is natural for us to conjecture that it may be also
true for ($p$,1)-total labellings. That is,
$C_{p,1}^T(G)=\chi_{p,1}^T(G) (=\lambda_p^T(G)+1)$. Unfortunately,
we could find counterexamples with $C_{p,1}^T(G)$ is strictly
greater than $\chi_{p,1}^T(G)$. Taking $P_k$ with $k>p$ as an
example, we have $\chi_{p,1}^T(G)\leq p+3$ by Theorem 2.2 but
$C_{p,1}^T(P_k)\geq2p$ by Theorem 2.4, which is strictly greater
than $\chi_{p,1}^T(P_k)$ when $p\geq4$.

\vspace{2mm}

\par Although we can not present a conjecture like LECC or LTCC, we
may conjecture an upper bound for $C_{p,1}^T(G)$ for any graph $G$:
\begin{conj}
Let $G$ be a simple graph with maximum degree $\Delta$. Then
$$C_{p,1}^T(G)\leq \Delta+2p.$$
\end{conj}
\par Obviously, the conjecture is true for paths and trees by
Theorem 2.4 and 2.6. Havet and Yu [6] gave a similar conjecture on
$\lambda_p^T(G)$. They also showed that $\lambda_p^T(K_n)=n+2p-2$
for complete graph with $n\geq6p^2-10p+4$ was even. Then
$C_{p,1}^T(K_n)\geq\chi_{p,1}^T(K_n)=\lambda_p^T(K_n)+1=\Delta+2p$.
Therefore, the bound in Conjecture 2.8 is tight.

\bigskip

\section{\bf Stars}

\par In this section, we prove that the conjecture above is true for
stars. Actually, we can improve the bound by one for stars.
\par Obviously, $C_{1,1}^T(K_{1,n})=\chi_l''(K_{1,n})=n+1$.
When $n\leq2$, $K_{1,n}$ is equivalent to $P_{n+1}$, which condition
have been shown in Theorem 2.4. Therefore, we only need to consider
the case when $p\geq2$ and $n\geq3$.

\begin{center}
\includegraphics[height=3cm]{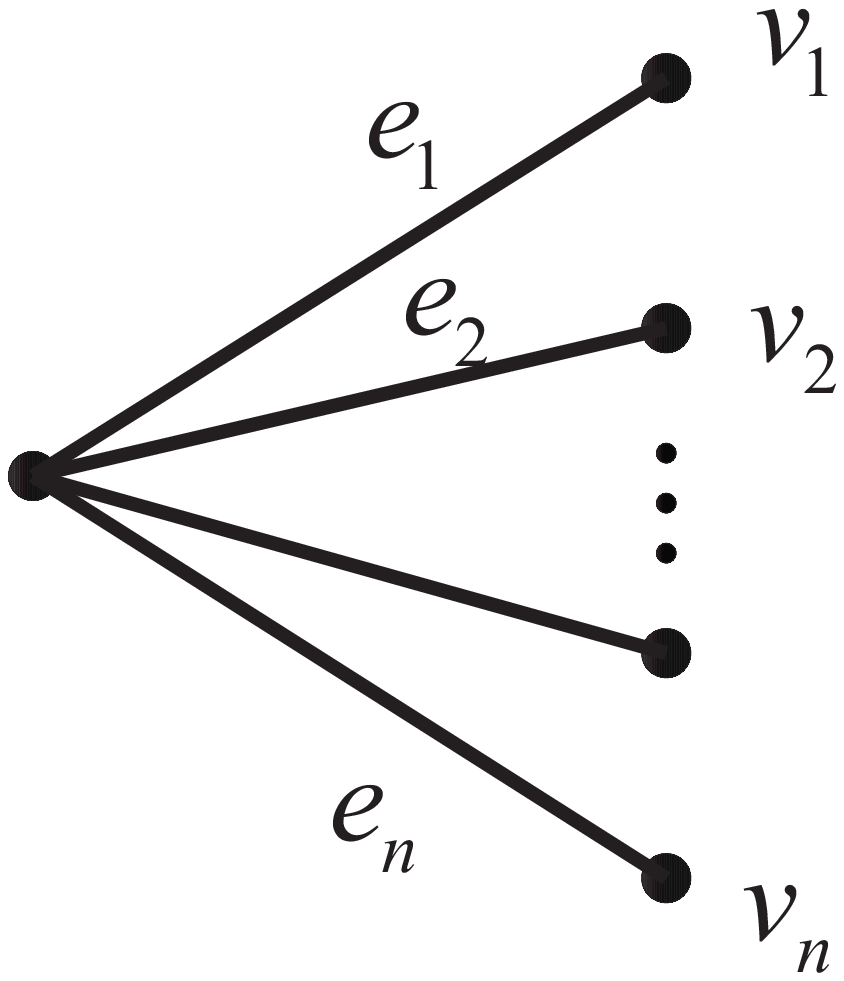}\\
\mbox{Figure 1}
\end{center}

\begin{thm}
Let $K_{1,n}$ be a star with $n\geq3$ and $p\geq2$. Then
$$C_{p,1}^T(K_{1,n})\leq n+2p-1.$$
\end{thm}

\noindent{\bf Proof.} Assume $|L(x)|=k$ with $k=n+2p-1$ for all
$x\in V\cup E$. Denote the maximum vertex by $w$ and the others by
$v_1, \cdots, v_n$. Denote the edges by $e_1, \cdots, e_n$,
respectively (see Figure 1). Denote the colors
$\{x-(p-1),\cdots,x-1,x,x+1,\cdots,x+(p-1)\}$ by $\|x\|_p$ and the
labelling of $K_{1,n}$ by $c$. Then if we label $x\in V\cup E$ with
color $\alpha\in L(x)$, we sometimes denote that by $c(x)=\alpha$.
\par First, label $w$ by the minimum color $\alpha$ of its list and
let $L'(e_j)=L(e_j)\setminus
\{||\alpha||_p\}$,$L'(v_j)=L(v_j)\setminus \{\alpha,||c(e_j)||_p\}$
for all $j\in\{1,2,\cdots,n\}$. Then we have $|L'(e_j)|\geq n-1$ and
$|L'(v_j)|\geq k-(2p-1+1)=n-1\geq1$. Therefore, we just need to
consider the coloring, denoted by $c$, of edges $e_j$ for all $j$.
The coloring of vertices $v_j$ is obvious since $|L'(v_j)|\geq1$.
Then we get an $L$-($p$,1)-total labelling of $K_{1,n}$ with the
assignment $L$.
\par Assume that at least one of the lists, say $L'(e_1)$, still
contains at least $n$ colors. We give an algorithm for the edge
coloring as follows:
\vspace{2mm}
\par Step 1: Let $i=1,S=\emptyset,a_i=e_1$;
\vspace{2mm}
\par Step 2: Determine the minimum color $m$ of the union of the
lists of all uncolored edges. That is, $m=\min\{x\mid x\in
\bigcup\limits_{e_p}L'(e_p)\}$ where $e_p\in E\setminus S$;
\vspace{2mm}
\par Step 3: If $L'(e_1)$ contains $m$ and no other uncolored edges
has $m$ in its list, then let $e_i'=a_i$; otherwise, choose another
$e_k$ with $m\in L'(e_k)$ and let $e_i'=e_k$.
\vspace{2mm}
\par Step 4: Let $c(e_i')=m$, $S=S\cup \{e_i'\}$;
\vspace{2mm}
\par Step 5: If $i=n$, then stop; otherwise, delete $m$ from the
lists of uncolored edges, that is, let $L'(e_p)=L'(e_p)\setminus
\{m\}$ for all $e_p\in E\setminus S$;
\vspace{2mm}
\par Step 6: If $e_i'=a_i$, then $a_{i+1}=e_p$ where $|L'(e_p)|\geq n-i$, $e_p\in E\setminus S$;
 else $a_{i+1}=a_i$;
\vspace{2mm}
\par Step 7: $i=i+1$, turn Step 2.
\vspace{2mm}
\par We delete at most one color in every step. So if $e_1$ is the
last edge colored by our algorithm, then the coloring is possible
since the list of $e_1$ has at least one color left by assumption.
If $e_1$ is not the last edge, then the coloring of $e_1$ deletes no
color from any list of $E\setminus S$. Suppose $e_1$ get colored by
the $i$th loop for some $i$. Then we have deleted at most $i-1$
colors from the list of $e_p$ for all $e_p\in E\setminus S$, and we
can choose some $e_p$ as the new beginning of our algorithm since we
have $|L'(e_p)|\geq n-i$ at the beginning of the next loop.
\par Thus, each edge list $L'(e_j)$ has exactly $n-1$ colors. That
means $||\alpha||_p \subseteq L(e_j)$ for all $j$. If we could not
finish the coloring, then an analogous fact must hold for every
color $\beta\in L(w)$. Therefore,
$\{\alpha-(p-1),\cdots,\alpha-1\}\cup L(w)\subseteq L(e_j)$. So we
have $k=|L(e_j)|\geq |\{\alpha-(p-1),\cdots,\alpha-1\}\cup
L(w)|=p-1+k$, which is a contradiction.\qed

\begin{lem} {\upshape ([6])} Let $G$ be a bipartite graph. Then
$$\Delta+p-1\leq\lambda_p^T(G)\leq\Delta+p.$$ Moreover, if
$p\geq\Delta$ or $G$ is regular, then $\lambda_p^T(G)=\Delta+p$.
\end{lem}

\begin{thm}Let $K_{1,n}$ be a star. Then
$$\chi_{p,1}^T(K_{1,n})=\left\{
\begin{array}{ll}
n+p, & \hbox{$p<n$;} \\
n+p+1, & \hbox{$p\geq n$.}
\end{array}
\right.$$
\end{thm}

\noindent{\bf Proof.} By Lemma 3.2 and $\chi_{p,1}^T=\lambda_p^T+1$,
we have $n+p\leq\chi_{p,1}^T(K_{1,n})\leq n+p+1$ and
$\chi_{p,1}^T(K_{1,n})=n+p+1$ when $p\geq n$. If $p\leq n-1$, then
we give a ($p$,1)-total labelling of $K_{1,n}$ with colors
$\{1,2,\cdots,n+p\}$. Suppose $K_{1,n}$ is defined as Figure 1. We
color $w$ with $n+p$\ and color\ $e_j$\ with\ $j$\ for all
$j\in\{1,2,\cdots,n\}$. After that we color $v_j$ with $p+j$ \ for
$j=1,\cdots,n-1$ and color $v_n$ with color 1. Since $n-1\geq p$,
this coloring is a proper ($p$,1)-total labelling of $K_{1,n}$.
Therefore, $\chi_{p,1}^T(K_{1,n})=n+p$ when $p< n$.\qed

\vspace{5mm}

\par When $p=2$, we have $\chi_{2,1}^T(K_{1,n})=n+2$ by Theorem 3.3.
Then $C_{2,1}^T(K_{1,n})\geq\chi_{2,1}^T(K_{1,n})=n+2$. On the other
hand, we also have $C_{2,1}^T(K_{1,n})\leq n+2\times2-1=n+2$ by
Theorem 3.1. That is, $C_{2,1}^T(K_{1,n})=n+2=n+2\times p-1$.
Therefore, the upper bound of Theorem 3.1 is tight when $p=2$.

\bigskip

\section{\bf Outerplanar graphs}

\par In this section, we discuss the $C_{p,1}^T$ of outerplanar
graphs $G$. An outerplanar graph is a planar graph that can be drawn
on the Eucliden plane such that there exists a face $f$ with all
$v\in V(G)$ belong to $f$. For these special graphs, we give a
theorem as follows:
\vspace{3mm}
\begin{thm}Let $G$ be an outerplanar graph with maximum degree $\Delta\geq p+3$. Then $$C_{p,1}^T(G)\leq \Delta+2p-1.$$
\end{thm}


\par We will prove Theorem 4.1 by contradiction. Before that, we
need a configuration lemma as follows:
\vspace{3mm}
\begin{lem} {\upshape ([3])} Every outerplanar graph $G$ with $\delta(G)=2$ contains one of the following configurations (see Figure 2):
\begin{itemize}
\item[$(C1)$] two adjacent 2-vertices $u$ and $v$;
\item[$(C2)$] a 3-face [$uv_1v_2$] with $d(u)=2$ and $d(v_1)=3$;
\item[$(C3)$] two 3-face [$u_1v_1x$] and [$u_2v_2x$] such that
$d(x)=4$ and $d(u_1)=d(u_2)=2$.
\end{itemize}
\end{lem}

\begin{center}
\includegraphics[height=3cm]{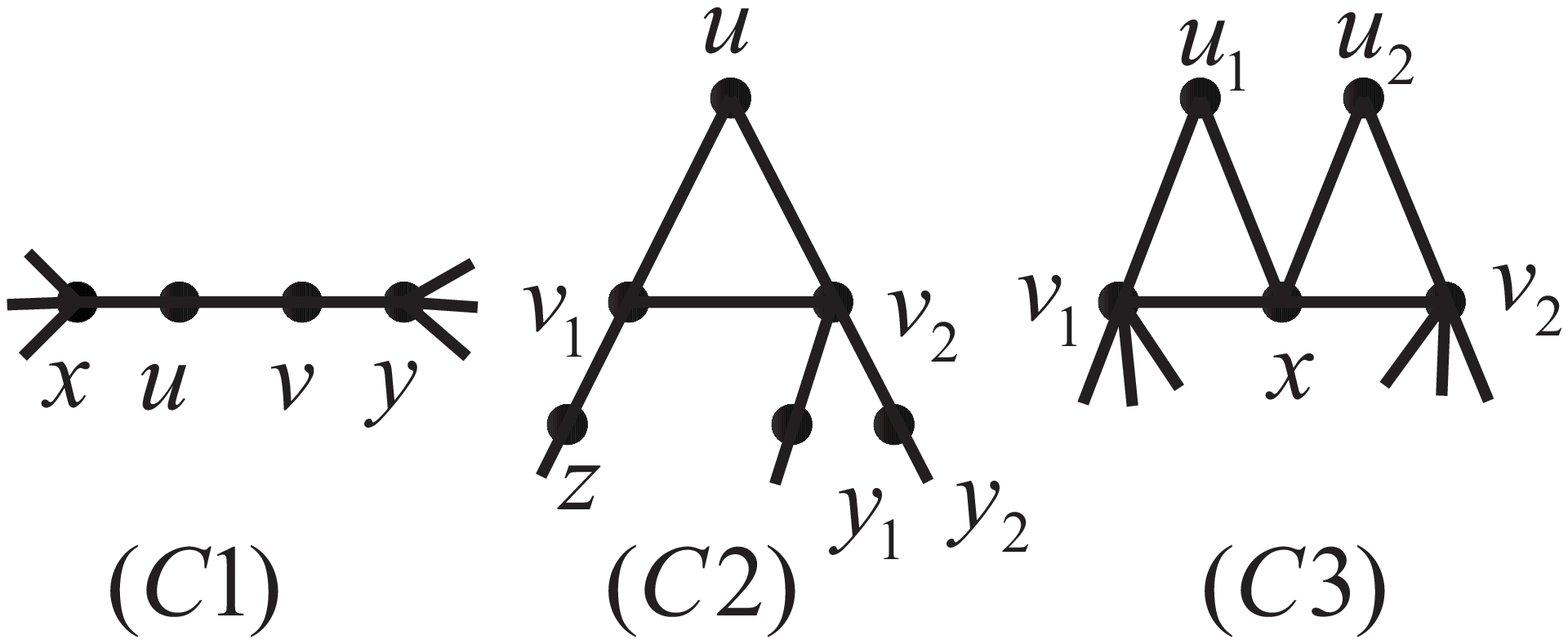}\\
\mbox{Figure 2}
\end{center}

\noindent{\bf Proof of Theorem 4.1.} Let $H$ be a minimal
counterexample in terms of $|V(G)|+|E(G)|$ to Theorem 4.1. $L$ is
the $k$ assignment defined on $V(H)\cup E(H)$ and $k=\Delta+2p-1$.
Denote the ($p$,1)-total labelling of $H$ by $c$. Then if we label
$x\in V(H)\cup E(H)$ with color $\alpha\in L(x)$, we sometimes
denote that by $c(x)=\alpha$. Denote by $L'(x)$ the set of colors
still available to color the element $x\in V(H)\cup E(H)$ such that
the labelling is a proper ($p$,1)-total labelling. We still use
$\|x\|_p$ to denote the color
set $\{x-(p-1),\cdots,x-1,x,x+1,\cdots,x+(p-1)\}$.\\

\vspace{2mm}

\noindent{\bf Claim 1.} $\delta(H)\geq 2$.\\

\noindent{\bf Proof.} If $\delta(H)=1$. Suppose that $e=uv\in E(H)$
and $d(v)=1$. The graph $H'=H\setminus \{v\}$ still satisfies the
demands of the theorem. By the minimality of $H$, $H'$ is
$L$-($p$,1)-total labelable. Without loss of generality, we suppose
the labelling is $c$. Then at most $\Delta-1+2p-1$ colors are
forbidden for the labelling of edge $e$. So we can choose a color
for $e$ from $L'(e)$ since $|L'(e)|\geq|L(e)|-(\Delta-1+2p-1)=1$.
After that we color $v$ from $L'(v)=L(v)\setminus
\{c(u),||c(e)||_p\}$. It is possible since $|L'(v)|\geq
k-(1+2p-1)=\Delta-1$. Then we extend the labelling $c$ to $H$, which
is a contradiction.\qed

\vspace{2mm}

Therefore, $\delta(H)=2$. By Lemma 4.2, $H$ contains one of the
configurations $C1$--$C3$. Next, we will show that in each case of
$C1$--$C3$ we can obtain a labelling such that $H$ is
$L$-($p$,1)-total labelable. Then we get a contradiction:\\

\noindent {\bf $(C1)$} Let $x$ be the neighbor of $u$ different from
$v$ and $y$ the neighbor of $v$ different from $u$. Let
$H'=H\setminus e$ where $e=uv$. Then $H'$ still satisfies the
demands of the theorem. By the minimality of $H$, $H'$ is
$L$-($p$,1)-total labelable. Remove the colors of vertex $u$ and
$v$. After that we define a list of available colors for each of
$u,v$ and $e$ as follows.
$$L'(u)=L(u)\setminus \{c(x),\|c(ux)\|_p\},$$
$$L'(v)=L(v)\setminus \{c(y),\|c(vy)\|_p\},$$
$$L'(e)=L(e)\setminus \{c(ux),c(vy)\}.$$
Since $|L|=k=\Delta+2p-1$ and $\Delta\geq p+3$, it follows that
$$|L'(u)|\geq k-(1+2p-1)\geq p+2,$$
$$|L'(v)|\geq k-(1+2p-1)\geq p+2,$$
$$|L'(e)|\geq k-2\geq 3p.$$
Let $m=\min \{x\mid x\in L'(u)\cup L'(v)\cup L'(e)\}$.\\
\noindent {\bf Case 1.}  $m\notin L'(u)\cup L'(v)$. Let $c(e)=m$ and
at most $p-1$ colors are unavailable for coloring $u,v$. Then at
least 3 colors are left in $L'(u)$ and $L'(v)$. So we can choose two
left colors from the list of $u,v$ such than $c(u)\neq c(v)$.\\
\noindent {\bf Case 2.} $m\in L'(u)$ or $L'(v)$. Without loss of
generality, say $m\in L'(u)$. Let $c(u)=m$. Then at most $p$ colors
are unavailable for coloring $e$ and at most one color for $v$. Let
$$L''(v)=L'(v)\setminus \{c(u)\},L''(e)=L'(e)\setminus
\{\|c(u)\|_p\}.$$ Then we have $$|L''(v)|\geq p+2-1\geq
p+1,|L''(e)|\geq 3p-(2p-1)\geq p+1.$$ Let $m_1=\min \{x\mid x\in
L''(v)\cup
L''(e)\}$.\\
\noindent {\bf Case 2.1.}  $m_1\in L''(v)$. Let $c(v)=m_1$ then we
delete at most $p$ colors from $L''(e)$. So we can color $e$ since
at least $|L''(e)|-p\geq 1$ colors are still available for $e$.\\
\noindent {\bf Case 2.2.} $m_1\notin L''(v)$. Let $c(e)=m_1$ then we
delete at most $p-1$ colors from $L''(v)$. So we can color $v$ since
at least $|L''(v)|-(p-1)\geq2$ colors are still available for $v$.\\
In any case, we extend the labelling $c$ to $H$ for $(C1)$, which is
a contradiction.\\

\noindent {\bf $(C2)$} Let $z$ be the neighbor of $v_1$ and $v_2$.
Let $H'=H\setminus uv_1$. Then $H'$ still satisfies the demands of
the theorem. By the minimality of $H$, $H'$ is $L$-($p$,1)-total
labelable. Remove the colors of vertex $u$. After that we define a
list of available colors for each of $u$ and $uv_1$ as follows.
$$L'(u)=L(u)\setminus \{c(v_1),c(v_2),\|c(uv_2)\|_p\},$$
$$L'(uv_1)=L(uv_1)\setminus \{c(v_1z),c(v_1v_2),c(uv_2),\|c(v_1)\|_p\}.$$
Since $|L|=k=\Delta+2p-1$ and $\Delta\geq p+3$, it follows that
$$|L'(u)|\geq k-(2+2p-1)\geq p+1,$$
$$|L'(uv_1)|\geq k-(3+2p-1)\geq p.$$
Let $m=\min \{x\mid x\in L'(u)\cup L'(uv_1)\}$. If $m\in L'(uv_1)$,
let $c(uv_1)=m$. Then at most $p$ colors are unavailable for
coloring $u$. So we can color $u$ since at least $|L'(u)|-p\geq1$
colors are still available for $u$; otherwise, let $c(u)=m$ and at
most $p-1$ colors are unavailable for coloring $uv_1$. So we can
color $uv_1$ since at least $|L'(uv_1)|-(p-1)\geq1$ colors are still
available for $uv_1$. Then we extend the labelling $c$ to $H$ for
$(C2)$, which is a contradiction.\\

\noindent {\bf $(C3)$} Let $H'=H\setminus xu_1$. By the minimality
of $H$, $H'$ is $L$-($p$,1)-total labelable. Remove the colors of
vertex $u_1$. After that we define the lists of available labels for
$u_1$ and $xu_1$ as follows.
$$L'(u_1)=L(u_1)\setminus \{c(v_1),c(x),\|c(u_1v_1)\|_p\},$$
$$L'(xu_1)=L(xu_1)\setminus \{c(v_1x),c(v_1u_1),c(xv_2),c(xu_2),\|c(x)\|_p\}.$$
Since $|L|=k=\Delta+2p-1$ and $\Delta\geq p+3$, it follows that
$$|L'(u_1)|\geq k-(2+2p-1)\geq p+1,$$
$$|L'(xu_1)|\geq k-(4+2p-1)\geq p-1.$$
If $|L'(xu_1)|\geq p$, we can color $xu_1$ and $u_1$ as we have done
in Case $(C2)$; otherwise, we have $|L'(xu_1)|=p-1$ and
$\Delta=p+3$. Let $$m=\min \{x\mid x\in L'(u_1)\cup L'(xu_1)\},$$
$$M=\max\{x\mid x\in L'(u_1)\cup L'(xu_1)\},$$
$$L_0'(u_1)=L'(u_1)\setminus \{m,M\}.$$
Denote by $m_1$ $(M_1)$ and $a$ $(b)$ the minimum (maximum) number
of $L_0'(u_1)$ and $L'(xu_1)$, respectively.\\

\par If $P=xu_1$ has not a partial list $L$-($p$,1)-total labelling for $u_1$ and $xu_1$, \ then we have some claims as follows.\\

\noindent{\bf Claim 2.} $L'(xu_1)$ is a series of $p-1$ successively
integers. That is, $|L'(xu_1)|=p-1$ and $b-a+1=p-1$.\\

\noindent{\bf Proof.} If $m\in L'(xu_1)$, let $c(xu_1)=m$. Then at
most colors $\{m,m+1,\cdots,m+(p-1)\}$ are forbidden for list
$L'(u_1)$. Therefore, we can color $u_1$ with at least one color
from the available colors of $L'(u_1)$. If $M\in L'(xu_1)$, then we
can finish the partial list-($p$,1)-total labelling for $u_1$ and
$xu_1$ with an analogous analysis. Therefore, we have $a\geq m+1$.
 If $b-a+1\geq p$, then $b-m\geq b-(a-1)\geq p$. We can color $u_1$ with $m$ and $xu_1$ with $b$.
Obviously, this is a partial list-($p$,1)-total labelling for $u_1$
and $xu_1$, which is a contradiction to our assumption. Therefore, $b-a+1=p-1$.\qed\\

\noindent{\bf Claim 3.} $L'(xu_1)=L_0'(u_1)$ and $m_1=m+1, M_1=M-1$.\\

\noindent{\bf Proof.}  If $a\leq m_1-1$ or $b\geq M_1+1$, then we
have $M-a\geq (M_1+1)-(m_1-1)=(M_1-m_1+1)+1\geq p$ or
$b-m\geq(M_1+1)-(m_1-1)=(M_1-m_1+1)+1\geq p$. We can color $u_1$
with $m$ and $xu_1$ with $b$ or we can color $u_1$ with $M$ and
$xu_1$ with $a$. Then we obtain a partial list-($p$,1)-total
labelling for $u_1$ and $xu_1$. Therefore, $a\geq m_1$ and $b\leq
M_1$. If $a\geq m_1+1$ or $b\leq M_1-1$, then we also have $b-m \geq
b-(m_1-1)\geq b-a+2=p$ or $M-a\geq (M_1+1)-a\geq(b+1)+1-a=p$ since
$b-a+1=p-1$ by Claim 2. We still obtain a partial list-($p$,1)-total
labelling for $u_1$ and $xu_1$. Therefore, $a=m_1$ and $b=M_1$. By
our assumption, $M_1-m_1+1\geq|L_0'(u_1)|= |L'(u_1)|-2\geq p-1$.
Thus, we have $M_1-m_1+1=|L_0'(u_1)|=p-1$. Together with $a=m_1,
b=M_1$ and Claim 2, we obtain $L'(xu_1)=L_0'(u_1)$.\qed\\

\noindent{\bf Claim 4.} $m_1=m+1$ and $M_1=M-1$.\\

\noindent{\bf Proof.} If $m_1\geq m+2$ or $M_1\leq M-2$, then we
have $b-m\geq b-(m_1-2)=M_1-m_1+2=p$ or
$M-a\geq(M_1+2)-a=M_1-m_1+2=p$ by Claim 3. We can color $u_1$ with
$m$ and $xu_1$ with $b$ or we can color $u_1$ with $M$ and $xu_1$
with $a$. Then we obtain a partial list-($p$,1)-total labelling for
$u_1$ and $xu_1$, which is a
contradiction to our assumption.\qed\\

That is to say,
$$L'(u_1)=\{m,m+1,\cdots,m+(p-1),m+p\},$$
$$L'(xu_1)=\{m+1,\cdots,m+(p-1)\}$$
where $a=m_1=m+1, b=M_1=m+(p-1), M=m+p$.
\par Claim 2--4 show that
if $L'(u_1), L'(xu_1)$ are not defined as above, we can obtain a
partial list-($p$,1)-total labelling for $u_1$ and $xu_1$. Next, we
show that we can obtain a partial list-($p$,1)-total labelling for
$u_1$ and $xu_1$ even if $L'(u_1), L'(xu_1)$ satisfies Claim 2--4 :
\vspace{3mm}
\par Since $|L'(xu_1)|=p-1$, it follows that $\Delta=p+3$ and
$\{c(v_1x),c(v_1u_1),c(xv_2),c(xu_2),\|c(x)\|_p\}$ are all distinct.
Otherwise, $|L'(xu_1)|>\Delta+2p-1-(2p+3)=p-1$ which is a
contradiction.
\par In particular,
$c(v_1u_1)\notin\{c(v_1x),c(xv_2),c(xu_2),\|c(x)\|_p\}$. We
interchange the colors of $xv_1$ and $u_1v_1$. After that, we define
a new list $L''$ of available colors for $u_1$ and $xu_1$, then
$$L''(xu_1)=L'(xu_1)=L(xu_1)\setminus \{c(v_1x),c(v_1u_1),c(xv_2),c(xu_2),\|c(x)\|_p\},$$
$$L''(u_1)=L(u_1)\setminus \{c(v_1),c(x),\|c(xv_1)\|_p\}.$$
\par Since $c(xv_1)\neq c(u_1v_1)$, we see that $L''(u_1)\neq
L'(u_1)$ and $|L''(u_1)|\geq |L'(u_1)|$. Then we have $m'\leq m-1$
or $M'\geq M+1$ where $m',M'$ denote the minimum and maximum number
of $L''(u_1)$. Otherwise, we have $m'\geq m$, $M'\leq M$ and
$|L''(u_1)|\leq M'-m'+1\leq M-m+1=|L'(u_1)|$. Since $|L''(u_1)|\geq
|L'(u_1)|$, $|L''(u_1)|=|L'(u_1)|$ and $m'= m$, $M'= M$. That is,
$L''(u_1)=L'(u_1)$,which is a contradiction.
\par If $m'\leq m-1$ then we have $m+(p-1)-m'\geq
m+(p-1)-(m-1)=p$. Let $c(u_1)=m'$ and $c(xu_1)=m+(p-1)$. Then we
obtain a partial list-($p$,1)-total labelling for $u_1$ and $xu_1$.
If $M'\geq M+1$, then we have $M'-a=M'-(m+1)\geq
(M+1)-(m+1)=(p+1)-1=p$. Let $c(u_1)=M'$ and $c(xu_1)=a=m+1$. Then we
obtain a partial list-($p$,1)-total labelling for $u_1$ and $xu_1$.
\par Any way, we extend the labelling $c$ to $H$ for $(C3)$, which is
a contradiction.\qed\\

\par Then we complete the proof of Theorem 4.1.\qed\\

\par When $p=2$, our result generalized a result of Chen and Wang [3]:
\begin{cor} {\upshape ([3]Theorem 7)} If $G$ is an outerplanar graph with $\Delta(G)\geq
5$, then $\lambda_2^T(G)\leq\Delta(G)+2$.
\end{cor}

\noindent{\bf Proof.} Obviously, when $\Delta\geq p+3$ we have
$\chi_{p,1}^T(G)\leq C_{p,1}^T(G)\leq \Delta+2p-1$ by Theorem 4.1.
Since $\chi_{p,1}^T(G)=\lambda_p^T(G)+1$, let $p=2$, we have
$\lambda_2^T(G)=\chi_{2,1}^T(G)-1\leq\Delta+3-1=\Delta+2$ with
$\Delta\geq5$.\qed

\vspace{2mm}

\par In [3], the author showed that there existed
infinitely many outerplanar graphs $G$ such that
$\chi_{2,1}^T(G)=\Delta+3$. So we have $C_{2,1}^T(G)=\Delta+3$ by
Theorem 4.1. That is to say, the upper bound in Theorem 4.1 can not
be improved when $p=2$.
\par Finally, we conjecture that Theorem 4.1 is also true when $\Delta\leq
p+2$.
\begin{conj}
Let $G$ be an outerplanar graph. Then $$C_{p,1}^T(G)\leq
\Delta+2p-1.$$
\end{conj}

\end{document}